\input amstex
\documentstyle{amsppt}
\magnification=\magstep1
\pageheight{9truein}
\pagewidth{6.5truein}

\NoBlackBoxes
\NoRunningHeads
%\NoPageNumbers
\nologo

\def\lra{\longrightarrow}
\def\op{\operatorname}
\def\br{\Bbb R}
\def\fg{\frak g}
\def\oim{\Omega_{inv}(M)}
\def\iva{i_V(\alpha)}
\def\hag{H^*_G(M; \br)}
\def\ham{H^*(M; \br)}

\topmatter
\title Canonical Equivariant Extensions Using Classical Hodge Theory 
\endtitle
\author Christopher Allday
\endauthor
\abstract 
Lin and Sjamaar have used symplectic Hodge theory to obtain canonical 
equivariant extensions for Hamiltonian actions on closed symplectic 
manifolds that have the strong Lefschetz property.  Here we obtain canonical 
equivariant extensions much more generally by means of classical Hodge theory.
\endabstract
\endtopmatter

\document

\baselineskip=2\baselineskip

\specialhead 1.  Introduction \endspecialhead

In \cite{L, S}, Lin and Sjamaar show how to use symplectic Hodge theory to 
obtain canonical equivariant extensions of closed forms in Hamiltonian 
actions of compact connected Lie groups on closed symplectic manifolds which 
have the strong Lefschetz property.  In this paper, we show how to do the 
same using classical Hodge theory.  This has the advantage of applying far 
more generally.  Our method makes use of Green's operator, but, as we shall 
show in \cite{A, P2}, it is often possible to make explicit calculations.

For non--abelian compact connected Lie groups, we use the small model, which 
is much simpler than the Cartan model and which has been shown to be chain 
homotopy equivalent to the Cartan model by Alekseev and Meinrenken 
(\cite{A, M}).  In the abelian case the two models are the same. The final 
section, however, considers the Cartan model.

I would like to thank A. Alekseev and E. Meinrenken for many helpful 
discussions. In particular, the last section is due to them.

\specialhead 2.  Notation and Terminology \endspecialhead

Throughout this paper $G$ will denote a compact connected Lie group with 
Lie algebra $\frak g$; and $M$ will be a closed, connected, orientable, 
smooth manifold.  $G$ will be acting on $M$; and $M$ will be given an invariant 
Riemannian metric.  All cohomology will have real coefficients.

\demo{Definition (2.1)}  Let $i : M \lra M_G$ be the inclusion of a fibre in 
the Borel construction bundle $M_G \lra BG$.  $M$ (or the action) is said to 
have a cohomology extension of the fibre (CEF) if
$$
i^* : H^*_G(M; \br) \lra H^*(M; \br)
$$
is surjective.
\enddemo

\demo{Remarks}  `CEF' may also stand for `cohomologically extendable from the 
fibre' or `cohomological extendability of the fibre'.  Often `TNHZ' 
(`totally non--homologous to zero') has been used for this condition.  It 
implies that $H^*_G(M; \br)$ is a free $H^*(BG; \br)$--module.  And, when 
$G$ is a torus, CEF implies that
$$
\varphi^* : H^*_G(M; \br) \lra H^*_G(M^G; \br)
$$
is injective, where $\varphi : M^G \lra M$ is the inclusion of the fixed 
point set, $M^G$.  The injectivity follows from the Localization Theorem of 
Borel, Hsiang and Quillen.  (See, e.g., \cite{A, P1}.)  And the injectivity of 
$\varphi^*$ is sometimes expressed by saying that $M$ is (cohomologically) 
equivariantly formal (which one might abbreviate CEF).

The purpose of this paper is to give a canonical section of $i^*$ in the CEF 
case.  The method is an easy application of classical Hodge theory.

Using the invariant Riemannian metric we define the Hodge star operator, 
$*$, on $\Omega(M)$, and then $d^*$, the Laplacian (or Laplace--Beltrami) 
operator $\Delta$, and Green's operator $G$.  (It should be clear from 
the context when $G$ is the Lie group and when $G$ is Green's operator.)  
Since the metric is invariant, $*$, $d^*$, $\Delta$ and $G$, like $d$, restrict 
to operators on the invariant forms $\Omega_{inv}(M) = \Omega(M)^G$.  Thus 
the usual Hodge Decomposition Theorem (see, e.g., \cite{Wa}, 6.8, or 
\cite{We}, Chapter IV, Theorem 5.2) applies to $\Omega_{inv}(M)$ without 
alteration.

The cohomology of the classifying space, $H^*(BG; \br)$, is a polynomial ring:  
$H^*(BG; \br) = \br[t_1, \dots, t_r]$, where each $t_j$ has positive even 
degree and $r$ is the rank of $G$.  We shall often denote this ring by 
$R_G$.  When $G$ is a torus, each $t_j$ has degree 2; and $R_G$ can be 
identified with the polynomial ring (symmetric algebra) on the dual of the 
Lie algebra of $G$:  that is, $R_G = S(\frak g^*)$.  More generally, $R_G = 
S(\frak g^*)^G$, the ring of invariants under the dual of the adjoint action.

The {\bf small model} for computing the equivariant cohomology, $H^*_G(M, 
\br)$, is $R_G \otimes \Omega_{inv}(M)$ with differential
$$
d_G = I \otimes d - \partial\,,
$$
where $\partial = \sum\limits^r_{j = 1}t_j \otimes i_j$; and, for each $j$, 
$1 \le j \le r$, $i_j$ is an operator on $\Omega_{inv}(M)$ of degree 
$-\deg(t_j) + 1$.  

To describe each $i_j$ in more detail, one considers $\wedge(\frak g)^G$, the 
subalgebra of the exterior algebra on $\fg$ fixed by the adjoint action.  
$\wedge(\fg)^G$ is $\wedge(\Cal P)$, the exterior algebra on the primitives, 
and it can be identified with $H_*(G; \br)$.  For a chosen basis $\{c_1, \dots, 
c_r\}$ of $\Cal P$, let $\{s_1, \dots, s_r\}$ be the dual basis of generators 
of $H^*(G; \br)$, which is also an exterior algebra.  Then $t_j$ corresponds 
to $s_j$ under transgression in the universal principal bundle $EG \lra BG$.  
In the formula for $\partial$, $i_j$ is the inner product by the multivector 
field induced on $M$ by $c_j$ via the exponential map and the group action 
in the usual way.  When $G$ is a torus, all $c_j$ have degree one, and the 
small model is the same as the Cartan model.  When $G$ is non--abelian, 
however, it is far from obvious that the small model correctly computes 
$H^*_G(M; \br)$:  that it does indeed do so is a theorem of Alekseev and 
Meinrenken (\cite{A, M}). Also see \cite{A, M} for more details of the
construction of the small model.
\enddemo

\demo{Definition (2.2)}  In the small model, $R_G \otimes \Omega_{inv}(M)$, 
let 
$$
P = \big(I \otimes d^*G\big)\partial\,,
$$
where $G$ here is Green's operator.
\enddemo

So $P$ is an operator of degree zero.  We shall usually abbreviate $I \otimes 
d$ and $I \otimes d^*G$ simply as $d$ and $d^*G$; and so $d_G = d - \partial$ 
and $P = d^*G\partial$.

\specialhead 3.  The Circle Case \endspecialhead

In this short section we give canonical equivariant extensions when $G = 
S^1$.  This case is simpler than the general case and nicely illustrates the 
method.  The small model is the same as the Cartan model in this case, 
namely, $R_G \otimes \Omega_{inv}(M)$; and $R_G = \br[t]$, where $\deg(t) = 
2$.  Also $\partial = t \otimes i_V$, where $V$ is the vector field coming 
from the circle action.  The case rests on the following lemma.

\proclaim{Lemma (3.1)}  Assume that $M$ has a CEF.  (See Definition (2.1).)  
Let $\alpha \in \Omega_{inv}(M)$, and suppose that $di_V(\alpha) = 0$.  Then 
$i_V(\alpha) = d\beta$ for some $\beta \in \oim$.
\endproclaim

\demo{Proof}  Since $di_V(\alpha) = 0$, $d_Gi_V(\alpha) = 0$.  Let $x = 
[i_V(\alpha)]_G$, the equivariant cohomology class of $i_V(\alpha)$.  So $x$ 
restricts to zero on the fixed point set, $M^G$, which is non--empty by CEF.  
So $x = 0$, again by CEF.  (See the Remarks following Definition (2.1).)  
Hence $i^*(x) = 0$, where $i : M \lra M_G$ is the inclusion of a fibre.  Thus 
$i_V(\alpha) = d\beta$ for some $\beta \in \oim$.\qquad$\square$
\enddemo

\proclaim{Theorem (3.2)}  Assume that $M$ has a CEF.  Let $\alpha \in \oim$ 
be a closed form (i.e., $d\alpha = 0$).  Let
$$
\widehat\alpha = (I - P)^{-1}(\alpha) = \alpha + P(\alpha) + P^2(\alpha) 
+ \cdots + P^n(\alpha) + \cdots
$$
($P^m(\alpha) = 0$ for all $m$ such that $2m > \deg(\alpha)$.)  Then 
$d_G\big(\widehat\alpha\big) = 0$.

Hence the map $\alpha \longmapsto \big[\widehat\alpha\big]_G$, restricted 
to harmonic forms, is a canonical cohomology extension of the fibre 
$H^*(M; \br) \lra H^*_G(M; \br)$.  
\endproclaim

\demo{Proof}  $d_G\alpha = -\partial\alpha = -ti_V(\alpha)$, where we have 
abbreviated $t \otimes i_V$ as $ti_V$.  And $d\iva = -i_Vd(\alpha) = 0$.  
Hence, by Lemma (3.1), $\iva$ is a boundary.  So $\iva = dd^*G\iva$.  
Now $d_G(\alpha + P(\alpha)) = -t\iva + dP(\alpha) - ti_VP(\alpha) = 
-ti_VP(\alpha)$.

Inductively, suppose that $d_G(\alpha + P(\alpha) + \cdots + P^j(\alpha)) = 
-ti_VP^j(\alpha)$.  Then $di_VP^j(\alpha) = -i_VdP^j(\alpha) = -i_V\partial 
P^{j - 1}(\alpha) = 0$.  So, again by Lemma (3.1), $i_VP^j(\alpha)$ is a 
boundary; and so $i_VP^j(\alpha) = dd^*Gi_VP^j(\alpha)$.  Hence $d_G(\alpha 
+ P(\alpha) + \cdots + P^j(\alpha) + P^{j + 1}(\alpha)) = -ti_VP^{j + 1}
(\alpha)$.\qquad$\square$
\enddemo

\demo{Example (3.3)} Suppose that $M$ is symplectic and that the action 
is Hamiltonian. Let $\omega \in \oim$ be the symplectic form, and let $\mu$
be the moment map. So $d\mu = i_V(\omega)$. Suppose, further, that $\mu$ 
has been chosen so that its average value is zero. Thus, in the Hodge
decomposition, the harmonic part of $\mu$ is zero; and so $\mu = d^*dG(\mu)
= d^*Gd(\mu)$. So $P(\omega) = t\mu$; and
$\widehat\omega = \omega + t\mu$,
the usual equivariant extension of $\omega$.
\enddemo

\specialhead 4.  The Small Model \endspecialhead

In this section $G$ is any compact connected Lie group.  Using the small 
model, the main result looks the same.

\proclaim{Theorem (4.1)}  Assume that $M$ has a CEF.  Let $\alpha \in \oim$ 
be a closed form (i.e., $d\alpha = 0$).  Let 
$$
\widehat\alpha = (I - P)^{-1}(\alpha) = \alpha + P(\alpha) + P^2(\alpha) + 
\cdots + P^n(\alpha) + \cdots
$$
Then $d_G\big(\widehat\alpha\big) = 0$.

Hence the map $\alpha \longmapsto \big[\widehat\alpha\big]_G$, restricted 
to harmonic forms, is a canonical cohomology extension of the fibre $H^*(M; 
\br) \lra H^*_G(M; \br)$.
\endproclaim

\demo{Proof}  The Localization Theorem (in a useful form) is not valid for 
non--abelian $G$.  We shall compensate for this by using induction not only 
on the power of $P$ but also on the degree of $\alpha$.  

First, $d_G(\alpha) = -\partial\alpha$.  By CEF, $\alpha$ has an equivariant 
extension; and so $\partial\alpha$ is a \break$d$--boundary.  Hence $\partial\alpha 
= dP(\alpha)$.  Since $\partial P(\alpha) = 0$ if $\deg(\alpha) \le 2$, 
this starts both inductions.  

Now suppose that $d_G(\alpha + P(\alpha) + \cdots + P^j(\alpha)) = 
-\partial P^j(\alpha)$.  For $I = (i_1, \dots, i_{j + 1})$, where 
$1 \le i_1 \le \cdots \le i_{j + 1} \le r = \op{rank}(G)$, let $t_I = t_{i_1} 
\cdots t_{i_{j + 1}}$.  Collecting terms we can set
$$
\partial P^j(\alpha) = \sum_It_I\alpha_I\,,
$$
where $\alpha_I \in \oim$ and $\deg(\alpha_I) = \deg(\alpha) - \deg(t_I) 
+ 1$.  So $\deg(\alpha_I) < \deg(\alpha)$.

Clearly $d\partial P^j(\alpha) = 0$.  So, since the monomials $t_I$ are 
linearly independent, $d\alpha_I = 0$ for all $I$.  Thus, by the induction 
hypothesis on degree, we can assume that $d_G\big(\widehat\alpha_I\big) = 0$ 
for all $I$.  On the other hand, for $i \ge 1$, 
$$
\sum_It_IP^i(\alpha_I) = P^i\sum_It_I\alpha_I = P^i\partial P^j(\alpha) = 0\,,
$$
because $P\partial = d^*G\partial^2 = 0$.  Thus 
$$
\sum_It_I\widehat\alpha_I = \sum_It_I\alpha_I = -d_G(\alpha + P(\alpha) + 
\cdots + P^j(\alpha))\,.
$$
Hence $\sum\limits_It_I\big[\widehat\alpha_I\big]_G = 0$.

Since $H^*_G(M; \br)$ is a free $R_G$--module by CEF, it follows that 
$i^*\big[\widehat\alpha_I\big]_G = 0$ for all $I$.  (As above, $i : M \lra 
M_G$ is the inclusion of a fibre.  And, for details, see Remarks (1) below.)  
Thus each $\alpha_I$ is a boundary; and so $\partial P^j(\alpha) = dP^{j + 1}
(\alpha)$.\qquad$\square$
\enddemo

\demo{Remarks}  (1)  The details of the free module argument are as follows.  
By CEF, there are $a_1, \dots, a_k$ in $H^*_G(M; \br)$ such that $\{a_1, 
\dots, a_k\}$ is a basis for $\hag$ as a free \break $R_G$--module, and $\{i^*(a_1), 
\dots, i^*(a_k)\}$ is a basis for $\ham$ as a vector space.  Let 
$\big[\widehat\alpha_I\big]_G = \sum\limits_j\lambda_{Ij}a_j + b_I$ where 
$\lambda_{Ij} \in \br$ and $b_I$ is a sum of terms of positive degree in 
$t_1, \dots t_r$.  Since $\sum\limits_It_I\big[\widehat a_I\big]_G = 0$, 
$\sum\limits_{I, j}
\lambda_{Ij}t_Ia_j = 0$, because all the $t_I$'s have the same polynomial 
degree (although not necessarily the same total degree).  Hence each 
$\lambda_{Ij} = 0$.  So $i^*\big[\widehat\alpha_I\big]_G = 0$ for each $I$.  
(In effect, for the purpose of this argument, we regrade $R_G$ so that 
$t_1, \dots, t_r$ all have the same degree, as in the torus case.)

(2) We shall give another proof of Theorem (4.1) in \cite{A, P2} using the 
minimal Hirsch--Brown model.

(3) In \cite{A, M}, Alekseev and Meinrenken give a canonical embedding of 
the small model into the Cartan model which is a homotopy equivalence of 
differential $R_G$--modules.  Applying this mapping to $\widehat\alpha$ gives 
a canonical equivariant extension of $\alpha$ in the Cartan model. In the 
next section, however, following suggestions by Alekseev and Meinrenken,
we obtain a version of Theorem (4.1) for the Cartan model directly.

\enddemo

\specialhead 5.  The Cartan Model\endspecialhead

Again, in this section, $G$ is any compact connected Lie group; but, here, 
we use the Cartan model instead of the small model.  Thus, in the definition 
of $P$, the operator $\partial$ is now that of the Cartan model.  The methods 
of this section, in particular Lemmas (5.3) and (5.4), are entirely due to 
A. Alekseev and E. Meinrenken.

Again, the theorem looks the same.

\proclaim{Theorem (5.1)}  Assume that $M$ has a CEF.  Let $\alpha \in 
\Omega_{inv}(M)$ be a closed form (i.e., $d\alpha = 0$).  Let
$$
\widehat\alpha = (I - P)^{-1}(\alpha) = \alpha + P(\alpha) + P^2(\alpha) 
+ \cdots + P^n(\alpha) + \cdots\ .
$$
Then $d_G\big(\widehat\alpha\big) = 0$.

Hence the map $\alpha \longmapsto \big[\widehat\alpha\big]_G$, restricted to 
harmonic forms, is a cohomology extension of the fibre $H^*(M; \br) \lra 
H^*_G(M; \br)$.
\endproclaim  

The theorem will follow directly from Lemma (5.4).  First, however, note that 
it is well--known that the inclusion of cochain complexes, $(\Omega_{inv}(M), 
d) \lra (\Omega(M), d)$, induces an isomorphism in cohomology.  It is only a little less well--known that the inclusion of cochain complexes
$$
\big(S(\fg^*)^G \otimes \Omega(M)^G, d\big) \lra 
\big((S(\fg^*) \otimes \Omega(M))^G, d\big)
$$
also induces an isomorphism in cohomology, $H^*(BG; \br) \otimes \ham$.  
($\Omega(M)^G = \oim$.)  Thus we have the next lemma.

\proclaim{Lemma (5.2)}  Let $a \in \big(S^p(\fg^*) \otimes \Omega^q(M)\big)^G$ 
and suppose that $da = 0$.  Then there is $b \in S^p(\fg^*)^G \otimes 
\Omega^q(M)^G$ and $c \in \big(S^p(\fg^*) \otimes \Omega^{q - 1}(M)\big)^G$ 
such that $a = b + dc$.  
\endproclaim

The next lemma shows that CEF implies the existence of more general 
equivariant extensions.  

\proclaim{Lemma (5.3)}  Assume that $M$ has a CEF.  Let $a \in \big(S^p(\fg^*) 
\otimes \Omega^q(M)\big)^G$, and suppose that $da = 0$.  Then there are 
$a_j \in \big(S^{p + j}(\fg^*) \otimes \Omega^{q - 2j}(M)\big)^G$ for 
$j \ge 0$, such that $a_0 = a$ and $da_j = \partial a_{j - 1}$ for all 
$j \ge 1$.  Thus $d_G(a + a_1 + \cdots + a_j + \cdots) = 0$.
\endproclaim

\demo{Proof}  Let $a = b + dc$ as in Lemma (5.2).  Then CEF clearly implies 
the existence of $b_j \in \big(S^{p + j}(\fg^*) \otimes \Omega^{q - 2j}(M)
\big)^G$ for $j \ge 0$, such that $b_0 = b$ and $db_j = \partial b_{j - 1}$ 
for all $j \ge 1$.  Now put $a_1 = b_1 - \partial c$, and put $a_j = b_j$ 
for $j \ge 2$.\qquad$\square$  
\enddemo

The final lemma, which easily implies Theorem (5.1), shows that arbitrary 
partial equivariant extensions can always be extended indefinitely 
(assuming CEF).

\proclaim{Lemma (5.4)}  Assume that $M$ has a CEF.  Let $a \in 
\big(S^p(\fg^*) \otimes \Omega^q(M)\big)^G$, and suppose that $da = 0$.  
Suppose given, for $0 \le j \le m$, $a_j \in \big(S^{p + j}(\fg^*) \otimes 
\Omega^{q - 2j}(M)\big)^G$ such that $a_0 = a$ and $da_j = \partial a_{j - 1}$ 
for $1 \le j \le m$.  Then there is $a_{m + 1} \in \big(S^{p + m + 1}(\fg^*) 
\otimes \Omega^{q - 2m - 2}(M)\big)^G$ such that $da_{m + 1} = \partial a_m$.  

In particular, one could take $a_{m + 1} = P(a_m)$.
\endproclaim

\demo{Proof}  We use induction on $m$.  The case $m = 0$ is clear by 
Lemma (5.3).  Also by Lemma (5.3), for $j \ge 1$, there are $y_j \in 
\big(S^{p + j}(\fg^*) \otimes \Omega^{q - 2j}\big)^G$ such that $dy_1 = 
\partial a$, and $dy_j = \partial y_{j - 1}$ for all $j \ge 2$.  For $1 \le 
j \le m$, let $c_j = a_j - y_j$.  Then $dc_1 = 0$, and for $2 \le j \le 
m$, $dc_j = \partial c_{j - 1}$.  Thus, by the induction hypothesis, there 
is $c_{m + 1}$ such that $dc_{m + 1} = \partial c_m$.  Now put $a_{m + 1} = 
c_{m + 1} + y_{m + 1}$.\qquad$\square$
\enddemo

\demo{Remark}  The analogue of Lemma (5.4) for the small model follows 
directly from the proof of Theorem (4.1).  And, thanks to Lemma (5.2), 
Theorem (5.1) could be proved along lines similar to those used to prove 
Theorem (4.1).  Equally, of course, Theorem (4.1) could be proven by the 
methods of this section.
\enddemo

\enddocument